\documentclass[12pt]{amsart}
\usepackage{amscd}

\newtheorem{thm}{Theorem}[section]
\newtheorem{cor}[thm]{Corollary}

\newtheorem{prop}[thm]{Proposition}
\theoremstyle{definition}

\newtheorem{rem}[thm]{Remark}
\newtheorem{rems}[thm]{Remarks}

\theoremstyle{remark}

\numberwithin{equation}{section}
\newcommand{\norm}[1]{\left\lVert#1\right\rVert}

\newcommand{\SL}{\mathit{SL}}
\newcommand{\R}{\mathbb{R}}
\newcommand{\C}{\mathbb{C}}
\newcommand{\Ha}{\mathbb{H}}
\newcommand{\trans}[1]{{}^t\kern-.2em{#1}}
\newcommand{\Trans}[1]{{}^T\kern-.2em{#1}}
\newcommand{\lsup}[2]{{}^{#1}\kern-.1em{#2}}
\newcommand{\simar}{\stackrel{\sim}{\rightarrow}}
\newcommand{\pair}[2]
{\pmb{\left\langle\right.}{#1},{#2}\pmb{\left.\right\rangle}}
\newcommand{\vol}{\mathit{vol}}

\newcommand{\Gm}{\Gamma}
\makeatletter
\def\simar{\ooalign{\hfil\raise.7ex\hbox{$\simeq$}\hfil
\crcr\raise.25ex\hbox{$\rightarrow$}}}
\makeatother
\makeatletter
\def\varin{\mathrel{\mathpalette\@varin\relax}}
\def\@varin#1{%
   \hbox{\setbox\z@\hbox{\m@th$#1\cup$}%
       \def\reserved@a{bold}%
       \dimen@\ifx\reserved@a\math@version .3\else .2\fi\p@
       \kern.5\wd\z@\kern-\dimen@
       \vrule\@width2\dimen@\@height1.08\ht\z@\@depth\z@
       \kern-\dimen@\kern-.5\wd\z@
       \box\z@}}

\DeclareFixedFont{\bgn}{OT1}{cmr}{m}{n}{20.74}
\DeclareFixedFont{\bgi}{OT1}{cmr}{m}{it}{20.74}
\newcommand{\bigzerol}{\smash{\hbox{\bgi O}}}
\newcommand{\bigzerou}{\smash{\lower1.7ex\hbox{\bgi O}}}

\DeclareMathOperator{\tr}{tr}

\def\bs#1{\boldsymbol{#1}}
\def\HLINE{\noalign{\hrule height1pt}}
\newcommand{\inner}[2]{\left({#1},{#2}\right)}
\newcommand{\Einner}[2]{\left\langle{#1},{#2}\right\rangle}
\makeatletter
\def\varin{\mathrel{\mathpalette\@varin\relax}}
\def\@varin#1{%
   \hbox{\setbox\z@\hbox{\m@th$#1\cup$}%
       \def\reserved@a{bold}%
       \dimen@\ifx\reserved@a\math@version .3\else .2\fi\p@
       \kern.5\wd\z@\kern-\dimen@
       \vrule\@width2\dimen@\@height1.08\ht\z@\@depth\z@
       \kern-\dimen@\kern-.5\wd\z@
       \box\z@}}
\def\eqnarray{%
   \stepcounter{equation}%
   \def\@currentlabel{\p@equation\theequation}%
   \global\@eqnswtrue
   \m@th
   \global\@eqcnt\z@
   \tabskip\@centering
   \let\\\@eqncr
   $$\everycr{}\halign to\displaywidth\bgroup
       \hskip\@centering$\displaystyle\tabskip\z@skip{##}$\@eqnsel
      &\global\@eqcnt\@ne \hfil$\displaystyle{{}##{}}$\hfil
      &\global\@eqcnt\tw@ $\displaystyle{##}$\hfil\tabskip\@centering
      &\global\@eqcnt\thr@@ \hb@xt@\z@\bgroup\hss##\egroup
        \tabskip\z@skip
      \cr}

\begin{document}
\title[Quantization of Geodesic Flow]
{Quantization of the Geodesic Flow \\ 
on Quaternion Projective Spaces} 



\author{Kenro Furutani}
\address{Kenro Furutani \endgraf
Department of Mathematics \endgraf 
Science University of Tokyo \endgraf 
2641 Noda, Chiba, Japan \endgraf 
postal code : 278-8510}
\email{furutani@ma.noda.sut.ac.jp}

\maketitle
\tableofcontents
\begin{abstract}
We study a problem of the geometric quantization for the quaternion
projective space.  First we explain a K\"ahler 
structure on the punctured cotangent bundle
of the quaternion projective space, whose K\"ahler form coincides with
the natural symplectic form on the cotangent bundle and show that
the canonical line bundle of this complex structure is holomorphically
trivial by explicitly constructing a nowhere vanishing holomorphic
global section.  Then we construct a Hilbert space consisting of a
certain class of holomorphic functions on the punctured cotangent
bundle by the method of
pairing polarization and incidentally we construct an operator
from this Hilbert space to the $L_2$ space of the quaternion
projective space.  Also we construct a similar operator between these
two Hilbert spaces through the Hopf fiberation.
We prove that these operators quantize
the geodesic flow of the quaternion projective space to
the one parameter group of the unitary Fourier integral operators
generated by the square root of the Laplacian plus suitable constant.
Finally we remark that the Hilbert space above has the reproducing kernel. 
\end{abstract}


\section*{Introduction}

In the paper ~\cite{FT} we explicitly constructed 
a K\"{a}hler structure
on the punctured cotangent bundles of complex and quaternion
projective spaces whose K\"ahler form coincides with the natural
symplectic form (see \cite{Ii2}, \cite{Sz1}, \cite{Sz2} and also a quite
recent paper \cite{Sz3}). 
This K\"{a}hler structure is also invariant under 
the action of the geodesic flow with respect to the Fubini-Study 
metric for the case of complex projective spaces and the metric 
defined through the Hopf fibration
for the case of quaternion projective spaces, respectively. Then 
in the paper ~\cite{FY} we applied this structure to construct 
quantization operators of the geodesic 
flow on complex projective spaces, by pairing
polarizations (see ~\cite{Ra1}, ~\cite{Ii1}).  
In this paper we construct a similar
operator for the quaternion projective space by the same
method as \cite{Ra2} and \cite{FY} (Theorem \ref{th:3.2}). 
We also construct a quantization operator 
by making use of the Hopf fiberation 
$S^{4n+3} \rightarrow P^n\mathbb{H}$.  These two operators do not
coincide.

Most parts of this paper are devoted to the concrete 
determination of fiber integrations in terms of the Gamma function. 
These arise from the pairing polarization (K\"ahler polarization 
and the vertical polarization) on the punctured cotangent
bundle of the quaternion projective space.

Similar to the cases of the sphere and the complex projective space, 
the operators we construct here are not unitary, however, in a sense, 
asymptotically unitary 
(modulo a constant multiple (= $\frac{\sqrt{2}}{\pi}$, 
independent of the dimension, Proposition \ref{prop:3.8})).

The case of the sphere was treated earlier by \cite{Ra2}.  
What remains among the compact symmetric spaces of rank one for which 
we want to construct a quantization operator of the geodesic flow, 
or more precise, an exact quantization operator of the
bicharacteristic flow of the square root of the Laplacian is only the 
case of the Cayley projective plane. 
We will treat that case in a separate paper. See also \cite{F}.

In $\S\,1$ we summarize standard properties of the quaternion
projective space and introduce a K\"ahler structure on the 
punctured cotangent bundle of the quaternion projective space 
with a slight modification from  \cite{FT}. In $\S\, 2$
we give explicit calculations of the pairing of the K\"ahler
polarization and the vertical polarization on the punctured 
cotangent bundle of the quaternion projective space. Especially 
we give a relation arising from the Hopf fibration between the Liouville
volume forms on the cotangent bundle of the sphere and the quaternion
projective space.

Based on these data, in $\S\, 3$ we construct 
a quantization operator which maps a certain
class of classical observables to $L_2$-functions on the quaternion
projective space.
In $\S\, 4$ we construct another quantization operator
by making use of the Hopf fibration and give a relation between these two
quantization operators.
Finally in $\S\, 5$  we note the existence of the reproducing kernel
of the Hilbert space consisting of a certain class of 
holomorphic functions which is constructed in $\S\, 3$.  Here we can
represent it only in the form of power series.


\section{Quaternion projective space and a K\"ahler structure}

In this section we describe a K\"{a}hler structure on the punctured
cotangent bundle of the quaternion projective space. 

Let $\mathbb{H}$ be the quaternion number field over the real number
field $\R$,
which is generated by $\{\bs{e}_0,\bs{e}_1,\bs{e}_2,\bs{e}_3\}$ with the
relations $\bs{e}_i \bs{e}_j$ given by the table
\begin{equation}
\renewcommand{\doublerulesep}{.5pt}
\begin{array}{|c||c|c|c|c|}
\hline
 & \bs{e}_0 & \bs{e}_1 & \bs{e}_2 & \bs{e}_3\\
\HLINE
\bs{e}_0 & \bs{e}_0 & \bs{e}_1 & \bs{e}_2 & \bs{e}_3\\
\hline
\bs{e}_1 & \bs{e}_1 & -\bs{e}_0 & \bs{e}_3 & -\bs{e}_2\\
\hline
\bs{e}_2 & \bs{e}_2 & -\bs{e}_3 & -\bs{e}_0 & \bs{e}_1\\
\hline
\bs{e}_3 & \bs{e}_3 & \bs{e}_2 & -\bs{e}_1 & -\bs{e}_0\\
\hline
\end{array}
\end{equation}

In this paper we regard $\mathbb{H}^n$ as a right $\mathbb{H}$-vector
space with the $\Ha$-inner product
\begin{equation}
\inner{h}{k}_{\mathbb{H}}=\sum_{i=1}^n\theta(h_i)k_i
\end{equation}
where $h=(h_1,\dots,h_n),k=(k_1,\dots,k_n)\in\mathbb{H}^n$ and
$\theta(x)=x_0\bs{e}_0-x_1\bs{e}_1-x_2\bs{e}_2-x_3\bs{e}_3$ for
$x = \sum\limits_{i=0}^3 x_i\bs{e}_i\in\mathbb{H}$. Then the $\R$
bilinear form 
$\Einner{h}{k}=\dfrac 12 \{\inner{h}{k}_{\mathbb{H}}
+\inner{k}{h}_{\mathbb{H}}\}$ on $\mathbb{H}^n$ 
defines a Euclidean inner product on
$\mathbb{H}^n$ as a real vector space.  We will denote its
extension to the complexification $\Ha^n\otimes\C$ as a complex bilinear
form with the same notation $\Einner{\,\cdot\,}{\,\cdot\,}$.

Let $M(n,\Ha)$ be the space of $n\times n$ $\mathbb{H}$-matrices, and
for $X=(x_{ij})\in M(n,\Ha)$ define respectively
\begin{align}
\theta(X) &=(\theta(x_{ij})),\\
(\trans{X})_{ij} &= x_{ji},\\
\tr X &= \sum x_{ii},\\
\Trans{X} &= \theta(\trans{X}).
\end{align}

Each matrix $X\in M(n,\Ha)$ defines a 
right $\Ha$-linear map
$X:\;\Ha^n\rightarrow\Ha^n$, and we have
\begin{equation}
\inner{Xh}{k}_{\Ha}=\inner{h}{\Trans{X}k}_{\Ha}.
\end{equation}
The group $Sp(n)$ is then defined as a group consisting of those matrices 
$X\in M(n,\Ha)$ which
preserve  the $\Ha$-inner product.
 
The space $\mathcal{H}(n,\Ha)$ = $\{X \in M(n,\Ha) | X = \Trans{X} \}$ 
is called a Jordan algebra with the
Jordan product
\begin{equation}
X\circ Y=\frac 12 (XY+YX).
\end{equation}

Note that for $X\in\mathcal{H}(n,\Ha)$, $\tr X\in\R \bs{e}_0=\R$
and $\mathcal{H}(n,\Ha)$ is equipped with a Euclidean inner product given by
\begin{equation}
\Einner{X}{Y}_{\R}=\tr(X\circ Y).
\end{equation}

The inner product has the property:
\begin{equation}\label{symmetric}
\Einner{X\circ Y}{Z}_{\R}=\Einner{X}{Y\circ Z}_{\R}.
\end{equation}

As is well-known, the quaternion projective space $P^n\Ha$ is the
set of all $\Ha$-one-dimensional subspaces in $\Ha^{n+1}$. Here we
identify $P^n\Ha$ with the subset in $\mathcal{H}(n+1,\Ha)$:
\begin{multline}
P^n(\Ha)=\{P\in\mathcal{H}(n+1,\Ha)\;|\; P=(p_i\theta(p_j)),\\
p=(p_0,\dots,p_n)\in \Ha^{n+1},\Einner{p}{p}=1\}.
\end{multline}

We consider the Riemannian metric on $P^n\Ha$ defined through the
Hopf-fibration $\pi:\;S^{4n+3}\rightarrow P^n\Ha$, where
$S^{4n+3}=\{h\in\Ha^{n+1}\;|\;\Einner{h}{h}=\inner{h}{h}_{\Ha}=1\}$
has the standard metric.

Let us denote the isomorphism $\Ha\rightarrow M(2,\C)$
\begin{equation}
\Ha\ni h\longmapsto\left(\begin{array}{@{\,}cc@{\,}}
x_0+\sqrt{-1}x_1 & x_2+\sqrt{-1}x_3\\
-x_2+\sqrt{-1}x_3 & x_0-\sqrt{-1}x_1
\end{array}\right)\in M(2,\C)
\end{equation}
by $\rho$ and we denote with the same notation $\rho$ its
complexification
\begin{equation}
\Ha\otimes\C\stackrel{\rho}{\simar} M(2,\C).
\end{equation}
\smallskip

Now we introduce the following spaces:

\begin{multline*}
\mathbf{E}_S =\{(p,q)\in\Ha^{n+1}\times\Ha^{n+1}\;|\;
\Einner{p}{p}=1,\Einner{p}{q}=0,\\ 
q + p(q,p)_{\Ha} \ne 0 \}.
\end{multline*}

The space $\mathbf{E}_S$ is $Sp(n+1)$-invariant and also invariant under the
right action of $Sp(1)$.
\smallskip

\begin{multline*}
\mathbf{E}_S^0 =\{(p,q)\in\Ha^{n+1}\times\Ha^{n+1}\;|\;
\Einner{p}{p}=1,q\ne 0,(q,p)_{\Ha}= 0 \},
\end{multline*}
\smallskip

\begin{multline*}
\mathbf{E}_{\Ha}=\left\{(P,Q)\;|\; P,Q\in\mathcal{H}(n+1,\Ha),\tr P=1,
P\circ Q=\frac 12 Q,Q\ne 0\right\},
\end{multline*}
\smallskip

\begin{multline*}
\widetilde{\mathbf{E}_S} =\left\{(B_0,\dots,B_n)
\in M(2,\C)^{n+1} \;|\; z \wedge w \ne 0, \sum\det B_i=0\right\},
\end{multline*}
where $B_i$ $\in$ $M(2,\C)$ takes the form 
$B_i = \begin{pmatrix}z_{2i}& w_{2i}\\ z_{2i+1} & w_{2i+1}\end{pmatrix}$,
and $z = (z_0, \ldots , z_{2n+1})$, $w = (w_0, \ldots ,w_{2n+1})$.
\smallskip

\begin{multline*}
\widetilde{\mathbf{E}_S^0} =\left\{(B_0,\ldots,B_n)
\in \widetilde{\mathbf{E}_S} \;|\; \sum B_i B^*_i 
= - J\sum {\overline  B_i}^tB_iJ \right\},
\end{multline*}
where $J \,=\begin{pmatrix}0 & 1\\ -1 & 0\end{pmatrix}$ and 
the condition $\sum B_i B^*_i = - J\sum {\overline  B_i}^tB_iJ$
can be rewritten as 
$\sum z_{2i}{\overline z_{2i+1}}+ w_{2i}{\overline w_{2i+1}} =0 $.
\smallskip

\begin{multline*}
\widetilde{\mathbf{E}_{\Ha}}=\{A\in M(2n+2,\C)\;|\;\mathbb{J}\,A=\,\trans{A}\,
\mathbb{J},\,\mathrm{rank}\,A = 2,\,A^2=0\}.
\end{multline*}

Here we denote by $\mathbb{J}$
\begin{align*}
\mathbb{J}&=\left(\begin{array}{@{\,}cccc@{\,}}
J & & & \bigzerou\\
 & J & &\\
 & & \ddots &\\
\bigzerol & & & J
\end{array}\right)\in M(2n+2,\C).
\end{align*}
\smallskip

We will define the norm of a matrix $A = (a_{ij}) \in M(n,\mathbb{C})$
by $\norm{A} = \sqrt{\sum |a_{ij}|^2}$ = $\sqrt{\tr(A\,A^*)}$ 
and the norm of an element $B\in \widetilde{\mathbf{E}}_S$
by $\norm{B} = \sqrt{\sum \norm{B_i}^2}$.

\begin{rems}\label{sec1:1}
The inner product $\left<\,\cdot\,,\,\cdot\,\right>_{\R}$ on
$\mathcal{H}(n+1,\Ha)$ is the restriction of the $\R$-bilinear form
$\dfrac 12\tr(X\Trans{Y}+Y\Trans{X})$ on $M(n+1,\Ha)$ and can be
extended to $M(n+1,\Ha)\otimes\C$ $\cong$ $M(2n+2,\C)$ 
as a complex
bilinear form in a natural way. We will denote this bilinear form by
\begin{equation*}
\Einner{\,\cdot\,}{\,\cdot\,}_{\C}:
M(n+1,\Ha)\otimes\C \times M(n+1,\Ha)\otimes\C \rightarrow \C.
\end{equation*}
Then the Hermitian inner product on
$M(n+1,\Ha)\otimes\C \cong M(2n+2,\C)$ 
is given by $\Einner{A}{\overline{B}}_{\C}$ 
and the norm of $A \in M(2n+2,\C)$ $\cong$ $M(n+1,\Ha)\otimes\C$,
which we introduced above,
equals $\norm A$ =$\sqrt{2\Einner{A}{\overline{A}}_{\C}}$,
that is $\Einner{A}{\overline{A}}_{\C}$ = $1/2\cdot\tr(A\,A^*).$
\end{rems} 

Next we define the maps 
$\alpha,\beta,\pi_{\Ha},\pi_S,\tau_S$ and $\tau_{\Ha}$
among these spaces:

{\allowdisplaybreaks
\begin{align*}
&\begin{array}{@{\,}ccccc@{\,}}
\alpha: & \mathbf{E}_S & \rightarrow & \mathbf{E}_{\Ha}&\\
 & \varin &&\varin&\\
 & (p,q) &\longmapsto & (P,Q), & P=(p_i\theta(p_j)),Q=(p_i\theta(q_j)
+q_i\theta(p_j))
\end{array}\\[3ex]
&\begin{array}{@{\,}ccccc@{\,}}
\beta: & \widetilde{\mathbf{E}_S} & \rightarrow & \widetilde{\mathbf{E}_{\Ha}}&\\
 & \varin &&\varin&\\
 & (B_0,\dots,B_n) &\longmapsto & A=(A_{ij}), & A_{ij}=-B_iJ\trans{B}_j J
\end{array}\\[3ex]
&\begin{array}{@{\,}cccc@{\,}}
\pi_S: & \mathbf{E}_S & \rightarrow & S^{4n+3}\\
 & \varin &&\varin\\
 & (p,q) &\longmapsto & p
\end{array}\\[3ex]
&\begin{array}{@{\,}cccc@{\,}}
\pi_{\Ha}: & \mathbf{E}_{\Ha} & \rightarrow & P^n\Ha\\
 & \varin &&\varin\\
 & (P,Q) &\longmapsto & P
\end{array}\\[3ex]
&\begin{array}{@{\,}cccc@{\,}}
\tau_S: & \mathbf{E}_S & \rightarrow & \widetilde{\mathbf{E}_S}\\
 & \varin &&\varin\\
 & (p,q) &\longmapsto & (B_0,\dots,B_n),
\end{array}\\
& \phantom{\tau_S:\quad}B_i=\rho(\norm{q}p_i\otimes 1+q_i\otimes\sqrt{-1}),\quad
\norm{q}=\sqrt{\inner{q}{q}_{\Ha}}\\
&\begin{array}{@{\,}cccc@{\,}}
\tau_{\Ha}: & \mathbf{E}_{\Ha} & \rightarrow & \widetilde{\mathbf{E}_{\Ha}}\\
 & \varin &&\varin\\
 & (P,Q) &\longmapsto & A=(A_{ij}),
\end{array}\\
& \phantom{\tau_{\Ha}:\quad}A=\norm{Q}^2(\rho(P_{ij})) 
- \left(\rho(Q_{ij})\right)^2 +
\frac{\sqrt{-1}}{\sqrt{2}}\norm{Q}\left(\rho(Q_{ij})\right).
\end{align*}}
\smallskip

Of course $\mathbf{E}_{\Ha}$ is identified with the punctured tangent 
bundle of $P^n\Ha$.  The Riemannian metric $g_{\Ha}$ on $P^n{\Ha}$
defined through the Hopf fibration is given by
\begin{equation}\label{riemann}
g_{\Ha}( Q_1, Q_2) = \frac{1}{2}\,<Q_1,\,Q_2>_{\mathbb{R}} = 
\frac{1}{2}\,\tr(Q_1\circ Q_2)
\end{equation}
for $(P,Q_1)$, $(P,Q_2)$ $\in$ $\mathbf{E}_{\Ha}$.

\begin{rems}\label{sec1:2}
\begin{enumerate}
\item Let $(p,q) \in \mathbf{E}_S$.
Then $4\norm {q}^2$ = $\norm {B}^2$($B = \tau_S(p,q))$.

\item Let $(P,Q) \in \mathbf{E}_{\Ha}$, 
then $Q^3 = 1/2\norm{Q}^2 Q$, $2\norm {Q}^4$
= $\norm {A}^2$($A = \tau_{\Ha}(P,Q)$).

\item Assume that $(p,q)_{\Ha} = 0$, that is $(p,q) \in \mathbf{E}_S^0$,
then we can easily show that 
$2\norm{q}^2 = \norm{Q}^2$($\alpha(p,q) = (P,Q)$),
$Q^2 = 2\norm{q}^2 P + (q_i\theta(q_j))$ and $\beta^*(\norm{A})$ 
= $\frac{1}{\sqrt{2}}\norm{B}^2$($A = \tau_{\Ha}(P,Q),\, B = \tau_S (p,q)$).
Also we have $g_{\Ha}(Q,Q)$ = $\norm{q}^2$ for $(P,Q) = \alpha (p,q)$
with $(p,q) \in \mathbf{E}_S^0$.
\end{enumerate}
\end{rems} 

$\mathbf{E}_S$ is an open subspace of the tangent bundle of $S^{4n+3}$.
It consists of those vectors which are not parallel to the fiber of
the Hopf-fibration.

The map $\tau_S$(resp. $\tau_{\Ha}$) is an isomorphism between
the spaces $\mathbf{E}_S$ and $\widetilde{\mathbf{E}_S}$
(resp. $\mathbf{E}_{\Ha}$ and $\widetilde{\mathbf{E}_{\Ha}}$).  Also we have
$\tau_S(\mathbf{E}_S^0)$ = $\widetilde{\mathbf{E}_S^0}$ and

\begin{prop}\label{prop:1.3}
\begin{enumerate}
\item $\mathbf{E}_S^0$ is $Sp(1)$-invariant, 
and $\widetilde{\mathbf{E}_S^0}$ is $SU(2)$-invariant,
\item the following diagram is commutative:
\begin{equation}
\begin{CD}
S^{4n+3} @<{\pi_S}<< \mathbf{E}_S^0 @>{\tau_S}>> \widetilde{\mathbf{E}_S^0}\\
@VV{\pi}V  @VV{\alpha}V @VV{\beta}V\\
P^n\Ha @<<{\pi_{\Ha}}< \mathbf{E}_{\Ha} @>{\tau_{\Ha}}>>
\;\widetilde{\mathbf{E}_{\Ha}}.
\end{CD}
\end{equation}
\item the map $\tau_S$
commutes (on $\mathbf{E}_S$) with the action of 
$Sp(1)$ $\stackrel{\rho}{\simar}$ $SU(2)$.
\end{enumerate}
\end{prop}

\begin{rem}
The maps $\beta\circ\tau_S$ and $\tau_{\Ha}\circ\alpha$ do not
coincide on the whole space $\mathbf{E}_S$.
\end{rem}

By the Riemannian metric we identify the tangent bundle of 
the sphere $S^{4n+3}$, respectively the quaternion
projective space $P^n\Ha$ with their cotangent bundle.
Both spaces $\mathbf{E}_S$ and $\mathbf{E}_{\Ha}$ can be
seen as complex manifolds through the maps $\tau_S$ and $\tau_{\Ha}$. 
Then we have

\begin{prop}\label{sec1:kaehler}
Let $\omega_S$ and $\omega_{\Ha}$ be the symplectic forms on the
cotangent bundle of the sphere $S^{4n+3}$ and the quaternion
projective space, then
\begin{align}
\omega_S &= \tau_S^*\left(\sqrt{-1}\,\overline{\partial}\partial
\sqrt{\sum_{i=0}^n\norm{B_i}^2}\right),\\
\omega_{\Ha} &= \tau_{\Ha}^* \left(2^{\frac 14}\sqrt{-1}\,\overline{\partial}
\partial\sqrt{\norm{A}}\right).
\end{align}
\end{prop}

In fact we have more precise relations:
\begin{prop}
\begin{equation}\label{sec1:diff1}
\sqrt{-1}\tau_{S}^*(\partial\norm{B}-\overline{\partial}
\norm{B})=2\,\theta_S
\end{equation}
\begin{equation}\label{sec1:diff2}                    
\sqrt{-1}\tau_{\Ha}^*(\partial\sqrt{\norm{A}}-\overline{\partial}
\sqrt{\norm{A}})=2^{\frac 34}\,\theta_{P^n\Ha}
\end{equation}
where $\theta_S$( resp. $\theta_{P^n\Ha}$) is the canonical one-form on
the cotangent bundle $T^*S^{4n+3}$(resp. $T^*P^n\Ha$).
\end{prop}

\begin{proof}
We only show \eqref{sec1:diff2}.  The formula \eqref{sec1:diff1} is proved more directly.

By using the formulas in Remarks \ref{sec1:1}, \ref{sec1:2} and a
property \eqref{symmetric} of the
Euclidean inner product in the Jordan algebra  $\mathcal{H}(n+1,\Ha)$,
we can make the following calculations:
\begin{align*}
&\tau_{\Ha}\,^*(\partial\Einner{A}{\overline{A}}_{\C}^{1/4}
- \overline{\partial}\Einner{A}{\overline{A}}_{\C}^{1/4})\\
&= 4\frac{\sqrt{-1}}{\sqrt{2}}\Einner{\norm{Q}^2P - Q^2}{d(\norm{Q}Q)}_{\C}\\
&=-\frac{\sqrt{-1}}{\sqrt{2}}\Einner{Q}{dP}_{\C}.
\end{align*}
Then we have the formula \eqref{sec1:diff2} by adjusting constants
according to the definition of the
Riemannian metric \eqref{riemann} on $P^n\Ha$.
\end{proof}

When we consider the action of 
$\SL(2,\C)=\rho(\{r \in\Ha\otimes\C \,|\, r\,\theta(r) = 1\})$ on
$\widetilde{\mathbf{E}_S}$ from the right, we have
\begin{prop}\label{prop:1.7}
\begin{equation}
\beta:\;\widetilde{\mathbf{E}_S}\rightarrow\widetilde{\mathbf{E}_{\Ha}}
\end{equation}
is a principal fiber bundle with the structure group $\SL(2,\C)$.
\end{prop}

\section{Pairing of polarization}

We will denote by $\mathcal{G}$(resp. $\mathcal{G}_S$) the positive complex 
polarization on $\widetilde{\mathbf{E}_{\Ha}}$  
(resp. $\widetilde{\mathbf{E}_S}$)
defined by the K\"{a}hler structure (Proposition \ref{sec1:kaehler}), 
that is the sub-bundle consisting of tangents of type $(0,1)$, 
and by $K^{\mathcal{G}}_{\Ha}$ (resp. $K^{\mathcal{G}}_S$) its
canonical line bundle.
In this section we will describe a nowhere vanishing holomorphic
global section of the canonical line bundle $K^{\mathcal{G}}_{\Ha}$, and then 
we determine explicitly the pairing of polarizations on 
$\widetilde{\mathbf{E}_{\Ha}}$.

Let $Z$ be a vector field on 
$\underbrace{M(2,\C)\times\dots\times M(2,\C)}_{n+1}\setminus\{0\}$ 
defined by
\begin{equation}
Z=\frac{1}{\norm{B}^2}
\left(\sum_{i=0}^{2n+1}\overline{\frac{\partial D}{\partial z_i}}
\frac{\partial}{\partial z_i}
+\sum_{i=0}^{2n+1}\overline{\frac{\partial D}{\partial w_i}}
\frac{\partial}{\partial w_i}\right)
\end{equation}
where we denote $B=(B_0,\dots,B_n)\in M(2,\C)\times\dots\times
M(2,\C)$ and $B_i=\begin{pmatrix}z_{2i} & w_{2i}\\
z_{2i+1} & w_{2i+1}\end{pmatrix}$ as before and $D$ = $\sum \det B_i$

We define a $(4n+3)$-form $\sigma_S$ on $M(2,\C)\times\dots\times
M(2,\C)\setminus\{0\}$ by
\begin{equation}
\sigma_S= \frac{1}{(2\sqrt{-1})^{2n+2}}
{\bf i}_Z\left(dz_{0}\wedge\cdots\wedge dz_{2n+1}
\wedge dw_{0}\wedge\cdots\wedge dw_{2n+1}\right)
\end{equation}
where ${\bf i}_Z$ denotes the interior product with the vector field $Z$.

Since $d\,D(Z)\, \equiv\, 1$,we have
\begin{equation}
 (2\sqrt{-1})^{2n+2}\cdot dD\wedge\sigma_S = 
dz_{0}\wedge\cdots\wedge dz_{2n+1}\wedge dw_{0}\wedge\cdots\wedge dw_{2n+1}
\end{equation}
and the restriction of $\sigma_S$ to $\sum\det B_i=0$ is
holomorphic and nowhere vanishing.  We will denote the restriction of $\sigma_S$
to $\widetilde{\mathbf{E}_S}$ with the same notation.
\smallskip

Let $\begin{pmatrix}\sqrt{-1} & 0\\ 0 & -\sqrt{-1}\end{pmatrix},
\begin{pmatrix}0 & 1 \\ -1 & 0\end{pmatrix}$ and
$\begin{pmatrix}0 & \sqrt{-1}\\ \sqrt{-1} & 0\end{pmatrix}$
be elements in 
\smallskip

\noindent
$\mathfrak{su}(2)$ $\subset$ $\mathfrak{sl}(2,\C)$. 
We denote by $Y_1,Y_2$ and
$Y_3$ the vector fields on $\widetilde{\mathbf{E}_S}$
corresponding to these three elements respectively 
defined by the action of $\SL(2,\C)$.

Let us consider the $4n$-form $\sigma$ on $\widetilde{\mathbf{E}_S}$ given by
\begin{equation}
\sigma={\bf i}_{Y_3}\circ {\bf i}_{Y_2}\circ {\bf i}_{Y_1}(\sigma_S).
\end{equation}
Then $\sigma$ is a holomorphic $4n$-form, is invariant under the
action of the group $Sp(n+1)$ on $\widetilde{\mathbf{E}_S}$ from the
left and is $\SL(2,\C)$-invariant from the right. Also note that 
$det(Ad_g)$ = $1$ for any $g \in SL(2,\C)$.  So there exists 
a unique nowhere vanishing
holomorphic $4n$-form $\sigma_{\Ha}$ on $\widetilde{\mathbf{E}_{\Ha}}$ such that
\begin{equation}\label{eq:bss}
\beta^*(\sigma_{\Ha})=\sigma.
\end{equation}
This $4n$-form $\sigma_{\Ha}$ gives a holomorphic trivialization of the canonical
line bundle $K^{\mathcal{G}}_{\Ha}$ of
$\widetilde{\mathbf{E}_{\Ha}}$, and is $Sp(n+1)$-invariant.

Let $V_1,V_2$ and $V_3$ be three vector fields on the sphere $S^{4n+3}$
corresponding to the elements $\bs{e}_1,\bs{e}_2$ and $\bs{e}_3$
in $\mathfrak{sp}(1)$ $\subset$ $\Ha$ defined through the action of $Sp(1)$
from the right. Let $\bs{v}_S$ be the volume form of the Riemannian
metric on $S^{4n+3}$ and denote the volume form on $P^n\Ha$
by $\bs{v}_{\Ha}$.

We define three one-forms $\eta_1,\eta_2$ and $\eta_3$ on $S^{4n+3}$
in such a way that
\begin{align}
\eta_i(V_j)&=\delta_{ij}\\
\eta_i(V) &=0\quad\text{for any $V\in TS^{4n+3}$}
\end{align}
which is orthogonal to $V_j\,(j=1,2,3)$.

Now we have the following relations among these vector fields,
one-forms and volume elements:
\begin{align}
&\pi_* (\bs{v}_S)=2\pi^2\bs{v}_{\Ha},\label{eq:vS}\\ 
&2\pi^2\pi^*(\bs{v}_{\Ha})={\bf i}_{V_3}\circ {\bf i}_{V_2}\circ {\bf i}_{V_1}
(\bs{v}_S),\label{eq:vH}\\ 
&\eta_1\wedge\eta_2\wedge\eta_3\wedge {\bf i}_{V_3}\circ {\bf i}_{V_2}\circ
{\bf i}_{V_1}(\bs{v}_S)=\bs{v}_S.\label{eq:eta-iV}
\end{align}
Here $\pi_*$ means the fiber integration of the Hopf bundle 
$\pi : S^{4n+3}\rightarrow P^n\Ha$.

Let $\theta_i=(\pi_S\circ\tau_S^{-1})^*\eta_i$. We decompose
$\theta_i$ into
\begin{equation}
\theta_i=\theta_i'+\theta_i''
\end{equation}
with the holomorphic component $\theta_i'$ and the anti-holomorphic
component $\theta_i''$ in the complexified cotangent bundle
$T^{*}(\widetilde{\mathbf{E}_S})\otimes \C$.

Then, we have
\begin{prop}\label{prop:sigmaS}
\begin{equation}
\theta_1'\wedge\theta_2'\wedge\theta_3'\wedge ({\bf i}_{Y_3}
\circ {\bf i}_{Y_2}\circ
{\bf i}_{Y_1}(\sigma_S))=\det(\theta_i'(Y_j))\sigma_S.
\end{equation}
on $\widetilde{\mathbf{E}_S}$.
\end{prop}

Put
\begin{equation}\label{eq:2.13}
\sigma_S\wedge\overline{\sigma}_S=A_S\Omega_S
\end{equation}
and
\begin{equation}\label{eq:2.14}
\sigma_{\Ha}\wedge\overline{\sigma}_{\Ha}=A_{\Ha}\Omega_{\Ha}
\end{equation}
where
\begin{equation}\label{eq:2.15}
\Omega_S=\frac{(-1)^{(4n+3)(2n+1)}}{(4n+3)!}\omega_S^{\,\,4n+3}
= \frac{-1}{(4n+3)!}\omega_S^{\,\,4n+3}
\end{equation}
and
\begin{equation}\label{eq:2.16}
\Omega_{\Ha}=\frac{(-1)^{2n(4n-1)}}{(4n)!}\omega_{\Ha}^{\,\,4n}
=\frac{1}{(4n)!}\omega_{\Ha}^{\,\,4n}
\end{equation}
are the Liouville volume forms on $T^* S^{4n+3}$ and
$T^*P^n\Ha$ respectively.  

Then the invariance of $\sigma_S$ and $\Omega_S$ under the transitive 
action of the group $SO(4n+4)$ on the unit sphere in the tangent bundle
$TS^{4n+3}$ and the invariance of $\sigma_{\Ha}$ and $\Omega_{\Ha}$ 
under the transitive action of $Sp(n+1)$ on the unit sphere 
in the tangent bundle $TP^n\Ha$ give
\begin{prop}\label{prop:2.2}
\begin{align}
A_S & = \bs{a}_S \norm{B}^{4n+1}\label{eq:aS}\\ 
A_{\Ha}&= \bs{a}_{\Ha} \norm{A}^{2n+2}\label{eq:2.18}
\end{align}
with two constants $\bs{a}_S$ and $\bs{a}_{\Ha}$.
\end{prop}

Again by the same reasons as above, we can put
\begin{align}
&(\pi_S\circ\tau_S^{-1})^*(\bs{v}_S)\wedge\overline{\sigma}_S =B_S
\Omega_S, \label{eq:2.19}\\
&\qquad\quad B_S = \bs{b}_S \norm{B}^{-1}\notag
\end{align}
and
\begin{align}
&(\pi_{\Ha}\circ\tau_{\Ha}^{-1})^*(\bs{v}_{\Ha})\wedge
\overline{\sigma}_{\Ha}= B_{\Ha}\Omega_{\Ha},\label{eq:bH}\\ 
&\qquad\quad B_{\Ha} = \bs{b}_{\Ha}\norm{A}\notag,
\end{align}
with two constants $\bs{b}_S$ and $\bs{b}_{\Ha}$.
Then we have
\begin{prop}\label{prop:2.3}
\begin{equation}
\beta^*\left(\frac{A_{\Ha}}{B_{\Ha}}\right)=\det(\theta_i'(Y_j))
\frac{A_S}{B_S}
\end{equation}
on $\widetilde{\mathbf{E}_S^0}$.
\end{prop}

\begin{sloppypar}
\begin{proof}
By \eqref{eq:bss} and Proposition \ref{prop:sigmaS},
\begin{eqnarray*}
&& \sigma_S\wedge\overline{\sigma}_S=A_S\Omega_S\\
&=& \frac{1}{\det(\theta_i'(Y_j))\det(\theta_i''(Y_j))}
\theta_1'\wedge\theta_2'\wedge\theta_3'\wedge\theta_1''\wedge
\theta_2''\wedge\theta_3''\wedge\beta^*(\sigma_{\Ha})\wedge
\overline{\beta^*(\sigma_{\Ha})}\\
&=& \frac{\beta^*(A_{\Ha})}{\det(\theta_i'(Y_j))\det(\theta_i''(Y_j))}
\theta_1'\wedge\theta_2'\wedge\theta_3'\wedge\theta_1''\wedge
\theta_2''\wedge\theta_3''\wedge\beta^*(\Omega_{\Ha}).
\end{eqnarray*}
Then 
\begin{multline*}
\det(\theta_i'(Y_j))\det(\theta_i''(Y_j))\,A_S\,\mathbf{i}_{Y_3}
\circ\mathbf{i}_{Y_2}\circ\mathbf{i}_{Y_1}(\Omega_S)\\
=\beta^*(A_{\Ha})\,\mathbf{i}_{Y_3}\circ\mathbf{i}_{Y_2}\circ
\mathbf{i}_{Y_1}(\theta_1'\wedge\theta_2'\wedge\theta_3'\wedge
\theta_1''\wedge\theta_2''\wedge\theta_3'')\wedge\beta^*(\Omega_{\Ha})
\end{multline*}
on $\widetilde{\mathbf{E}_S}$, since $\beta_*(Y_i)=0$.

On $\widetilde{\mathbf{E}_S^0}$ 
we have by Proposition \ref{prop:1.3}, \eqref{eq:vS},
\eqref{eq:vH} and \eqref{eq:eta-iV}
\begin{eqnarray*}
&& \frac{\det(\theta_i'(Y_j))\det(\theta_i''(Y_j))}{B_S}\,
A_S\,\mathbf{i}_{Y_3}\circ\mathbf{i}_{Y_2}\circ\mathbf{i}_{Y_1}
\left((\pi_S\circ\tau_S^{-1})^*(\bs{v}_S)\wedge\overline{\sigma_S}\right)\\
&=& \frac{A_S}{B_S}\,\det(\theta_i'(Y_j))\,\det(\theta_i''(Y_j))\,
\mathbf{i}_{Y_3}\circ\mathbf{i}_{Y_2}\circ
\mathbf{i}_{Y_1}\left((\pi_S\circ\tau_S^{-1})^*
(\eta_1\wedge\eta_2\wedge\eta_3\right.\\
&&\hspace*{15em}
\left.\wedge\mathbf{i}_{\bs{V}_3}\circ
\mathbf{i}_{\bs{V}_2}\circ\mathbf{i}_{\bs{V}_1}(\bs{v}_S))
\wedge\overline{\sigma}_S\right)\\
&=&\frac{2\pi^2 A_S}{B_S}\,\det(\theta_i'(Y_j))\,\det(\theta_i''(Y_j))\\
&&\hspace*{5em}
\mathbf{i}_{Y_3}\circ\mathbf{i}_{Y_2}\circ
\mathbf{i}_{Y_1}\left((\pi_S\circ\tau_S^{-1})^*\left(\eta_1\wedge\eta_2\wedge
\eta_3\wedge\pi^*(\bs{v}_{\Ha})\right)\wedge\overline{\sigma}_S\right)\\
&=&\frac{2\pi^2 A_S}{B_S}\,\det(\theta_i'(Y_j))\,
\mathbf{i}_{Y_3}\circ\mathbf{i}_{Y_2}\circ
\mathbf{i}_{Y_1}\left(\theta_1\wedge\theta_2\wedge\theta_3
\wedge\beta^*\circ (\pi_{\Ha}\circ\tau_{\Ha}^{-1})^*(\bs{v}_{\Ha})\right.\\
&&\hspace*{16em}
\left.\wedge\theta_1''\wedge\theta_2''\wedge\theta_3''\wedge
\overline{\beta^*(\sigma_{\Ha})}\right).
\end{eqnarray*}
Here $\theta_i$ should be understood as restricted to
$\widetilde{\mathbf{E}_S^0}$. Now we have the following equality on
$\widetilde{\mathbf{E}_S^0}$:
\begin{eqnarray*}
&&\beta^*(A_{\Ha})\mathbf{i}_{Y_3}\circ\mathbf{i}_{Y_2}\circ
\mathbf{i}_{Y_1}(\theta_1'\wedge\theta_2'\wedge\theta_3'\wedge
\theta_1''\wedge\theta_2''\wedge\theta_3'')\wedge\beta^*(\Omega_{\Ha})\\
&=&\frac{2\pi^2 A_S}{B_S}\det(\theta_i'(Y_j))\mathbf{i}_{Y_3}
\circ\mathbf{i}_{Y_2}\circ
\mathbf{i}_{Y_1}(\theta_1'\wedge\theta_2'\wedge\theta_3'\wedge
\theta_1''\wedge\theta_2''\wedge\theta_3'')\wedge\beta^*(B_{\Ha}
\Omega_{\Ha}).
\end{eqnarray*}
Since $\mathbf{i}_{Y_3}\circ\mathbf{i}_{Y_2}\circ
\mathbf{i}_{Y_1}(\theta_1'\wedge\theta_2'\wedge\theta_3'\wedge
\theta_1''\wedge\theta_2''\wedge\theta_3'')\ne 0$, we have
\begin{equation}
2\pi^2\frac{A_S}{B_S}\det(\theta_i'(Y_j))=\beta^*\left(
\frac{A_{\Ha}}{B_{\Ha}}\right)
\end{equation}
on $\widetilde{\mathbf{E}_S^0}$.
\end{proof}
\end{sloppypar}

By \eqref{eq:aS} to \eqref{eq:bH} and $\norm B^2$ = $\sqrt{2}\beta^*(\norm A)$ 
on $\widetilde{\mathbf{E}_S^0}$ we have
\begin{cor}\label{cor:2.4}
\begin{equation}
2\pi^2\frac{\bs{a}_S}{\bs{b}_S}\det(\theta_i'(Y_j))
=\left(\frac{1}{\sqrt{2}}\right)^{2n+1}\frac{\bs{a}_{\Ha}}{\bs{b}_{\Ha}}.
\end{equation}
Hence $\det(\theta_i'(Y_j))$ must be constant on
$\widetilde{\mathbf{E}_S^0}$.   
\end{cor}

Now we list the concrete values of these constants 
$\bs{a}_S$, ${\bs{b}_S}$, ${\bs{a}_{\Ha}}$, ${\bs{b}_{\Ha}}$ 
and $\det(\theta_i'(Y_j))$.
\begin{prop}
\begin{align}
 \bs{a}_S &= -\sqrt{-1},\\
 \bs{b}_S &= \sqrt{-1},\\
 \bs{a}_{\Ha} &=2^{n-2},\\
\det(\theta_i'(Y_j))&=2^{-3}= \det(\theta_i'(Y_j'))
=\det(\theta_i''(Y_j''))\,\,\text{on}\,\, \widetilde{\mathbf{E}_S^0},\\
 \bs{b}_{\Ha}&=-\frac{1}{\sqrt{2}\pi^2}.
\end{align}
\end{prop}
Here $Y_i'$ (resp. $Y_i''$) is the holomorphic (resp. anti-holomorphic)
part of the vector field $Y_i$.
Note that the value $\det(\theta_i'(Y_j))$ is not constant on
all of $\widetilde{\mathbf{E}_S}$.

The calculations for the determination of these five constants are
so tedious that we do not write down the details here. They are
made by evaluating both sides of \eqref{eq:2.13} \eqref{eq:2.14} and
\eqref{eq:2.19} at a specific point.  Especially to determine the
constant $\bs{a}_{\Ha}$, we make use of the fiber bundle structure stated in
Proposition \ref{prop:1.7} to consider a suitable local coordinate neighborhood in
the space $\widetilde{\mathbf{E}_{\Ha}}$.  
Then we use the above Corollary \ref{cor:2.4} to determine the
constant $\bs{b}_{\Ha}$.   

Although it appeared in the proof of Proposition \ref{prop:2.3}
we emphasize the relation between the Liouville forms on $T^*S^{4n+3}$
and $T^*P^n\Ha$: 
\begin{prop}
\begin{multline}
\Omega_S\\
=\frac{\beta^*(A_{\Ha})}{A_S}
\frac{1}{|\det(\theta_i'(Y_j))|^2}
\theta_1'\wedge\theta_2'\wedge\theta_3'\wedge\theta_1''\wedge
\theta_2''\wedge\theta_3''\wedge\beta^*(\Omega_{\Ha}).
\end{multline}
\end{prop}

\section{Quantization operator I}

{}From now on we will omit the maps $\tau_{\Ha}$ and $\tau_S$ for 
the sake of simplicity.  

In this section we construct an operator from a Hilbert space
consisting of a certain class of holomorphic functions on 
$\widetilde{\mathbf{E}_{\Ha}}$ to $L_2(P^n\Ha)$.
 
The symplectic form $\omega_{\Ha}$ defines a complex line bundle $L$ on
$\widetilde{\mathbf{E}_{\Ha}}$ (although this is topologically trivial)
with the connection $\nabla$. Since $d\theta_{\Ha}=\omega_{\Ha}$
there is a trivialization of $L$ by a section $\bs{s}_{\Ha}$ such
that the connection $\nabla$ is given as
\begin{equation}
\nabla_X(\bs{s}_{\Ha})=2\pi\sqrt{-1}\Einner{\theta_{\Ha}}{X}\bs{s}_{\Ha},
\end{equation}
and since $\theta_{\Ha}$ is real we can introduce an inner product 
on $L$ such that
\begin{equation}
\inner{\bs{s}_{\Ha}}{\bs{s}_{\Ha}}_L\equiv 1.
\end{equation}

Again, by 
$\omega_{\Ha}=2^{\frac 14}\sqrt{-1}d(\partial\sqrt{\norm{A}})$, 
we have another trivialization of $L$
by a global section $\bs{t}_{\Ha}$ such that the connection is
expressed as
\begin{equation}
\nabla_X (\bs{t}_{\Ha})=2\pi\sqrt{-1}\Einner{2^{\frac 14}\sqrt{-1}\partial
\sqrt{\norm{A}}}{\,X}\bs{t}_{\Ha}.
\end{equation}
Since  $\bs{t}_{\Ha}$ = $\varphi\, \bs{s}_{\Ha}$ with a nowhere
vanishing function $\varphi$, this function $\varphi$ 
must satisfy the equation 
\begin{equation}
d\log\varphi=2\pi\sqrt{-1}(\theta_{\Ha}-2^{\frac 14}\sqrt{-1}\partial
\sqrt{\norm{A}}).
\end{equation}
{}From \eqref{sec1:diff2}, we can take a solution of this equation
\begin{equation}
\varphi=e^{-\sqrt[4]{2}\pi\sqrt{\norm{A}}}.
\end{equation}
We denote this solution by $\varphi_0$ and fix $\bs{t}_{\Ha}$ 
= $\varphi_0\, \bs{s}_{\Ha}$ henceforth.

Next we consider the canonical line bundle $K_{\Ha}^{\mathcal{G}}$ and its square
root $\sqrt{K_{\Ha}^{\mathcal{G}}}$. For the canonical line bundle 
$K_{\Ha}^{\mathcal{G}}$ we can introduce a connection 
$\lsup{\mathcal{G}}{\nabla}$ "along
the polarization" $\mathcal{G}$ in such a way that
\begin{equation}
\lsup{\mathcal{G}}{\nabla_X}(f)=\mathbf{i}_X\circ df,\,\text{for}\,\,
 X \in \Gamma (\mathcal{G}),\, f\in\Gamma (K_{\Ha}^{\mathcal{G}}).
\end{equation}
Note that 
$\mathbf{i}_X(f) = 0$ for $X \in \Gamma (\mathcal{G})$, 
$f\in \Gamma(K_{\Ha}^{\mathcal{G}})$.  
Now, by making use of this connection,
we can also define in a unique way a connection on the 
square root $\sqrt{K_{\Ha}^{\mathcal{G}}}$.
We denote this connection by 
$\lsup{\mathcal{G}}{\nabla^{\frac 12}}$.

{}From the pairing \eqref{eq:2.14} and by \eqref{eq:2.18} we 
define the pairing of 
$\sqrt{\sigma_{\Ha}}$ by
\begin{equation}
<\sqrt{\sigma_{\Ha}},\,\sqrt{\sigma_{\Ha}}> = 
\sqrt{|\bs{a}_{\Ha}|}\norm{A}^{n+1},
\end{equation}
and introduce an inner product for 
$\phi =f\cdot\bs{t}_{\Ha}\otimes\sqrt{\sigma_{\Ha}}$ $\in$
$\Gamma (L\otimes\sqrt{K_{\Ha}^{\mathcal{G}}})$
and 
$\psi = g\cdot\bs{t}_{\Ha}\otimes\sqrt{\sigma_{\Ha}}$ 
$\in$ $\Gamma (L\otimes\sqrt{K_{\Ha}^{\mathcal{G}}})$ 
($f$ and $g$ are functions on $\widetilde{\mathbf{E}_{\Ha}}$)
by
\begin{equation}\label{eq:3.9}
\inner{\phi}{\psi}
=\int_{\widetilde{\mathbf{E}_{\Ha}}} \,f(A)\overline{g(A)}\,
(\bs{t}_{\Ha},\,\bs{t}_{\Ha})_L<\sqrt{\sigma_{\Ha}},\,\sqrt{\sigma_{\Ha}}>\,
\Omega_{\Ha} 
\end{equation}
\begin{equation*}
=\int_{\widetilde{\mathbf{E}_{\Ha}}} \,f(A)\overline{g(A)}\,
e^{-2\,\sqrt[4]{2}\pi\sqrt{\norm{A}}}
|\bs{a}_{\Ha}|^{\frac 12}\,\norm{A}^{n+1}\Omega_{\Ha}.
\end{equation*}

Let us denote by
$\Gamma_{\mathcal{G}}(L\otimes\sqrt{K_{\Ha}^{\mathcal{G}}})$ 
the space of
smooth sections which are parallel with respect to the connection
$\nabla\otimes \mathit{Id} + \mathit{Id}\otimes
\lsup{\mathcal{G}}{\nabla^{\frac 12}}$. 
Then we have the isomorphism
\begin{equation}\label{eq:3.10}
f\longleftrightarrow f \cdot\bs{t}_{\Ha}\otimes\sqrt{\sigma_{\Ha}}
\end{equation}
between the space of holomorphic functions on $\widetilde{\mathbf{E}_{\Ha}}$
and $\Gamma_{\mathcal{G}}(L\otimes\sqrt{K_{\Ha}^{\mathcal{G}}})$. 
Note here that $\nabla_X(\bs{t}_{\Ha}) = 0$ for 
$X \in \Gamma(\mathcal{G})$.
We denote by
$\mathfrak{h}_{\Ha}^{\mathcal{G}}$ 
the completion of the space
\begin{align*}
&\{f\cdot\,\bs{t}_{\Ha}\otimes\sqrt{\sigma_{\Ha}}\,|\,f \,\text{is a
polynomial on}\, M(2n+2,\mathbb{C})\, \text{restricted to}\, 
\widetilde{\mathbf{E}_{\Ha}}\}\\
& \cong
\sum_{l=0}^{\infty}\mathcal{P}_l^{\Ha}
\end{align*}
with respect to the norm defined by the inner product 
\eqref{eq:3.9}, where $\mathcal{P}_l^{\Ha}$ denotes
the space of polynomials of degree $l$ on
$M(2n+2,\mathbb{C})$(restricted to $\widetilde{\mathbf{E}_{\Ha}}$). 
The spaces $\mathcal{P}_l^{\Ha}$ and
$\mathcal{P}^{\Ha}_{l^{'}}$ are orthogonal for $l \ne l^{'}$ with
respect to this inner product \eqref{eq:3.9}. 

Let $\mathcal{F}$ be the vertical polarization on
$\widetilde{\mathbf{E}_{\Ha}}\cong T_0^* P^n\Ha$ defined by the
projection 
$\pi_{\Ha}:\;T_0^*P^n\Ha\longrightarrow P^n\Ha$. In this case 
we also have a connection "along the polarization" $\mathcal{F}$ 
for the square root of the canonical line
bundle $K^{\mathcal{F}}_{\Ha}$ in the same way as for the case of 
the K\"{a}hler polarization
$\mathcal{G}$, and we have the correspondence
\begin{equation}\label{eq:3.11}
g\longleftrightarrow \pi_{\Ha}\,^*(g)\cdot\bs{s}_{\Ha}
\otimes\sqrt{\pi_{\Ha}\,^*(\bs{v}_{\Ha})}
\end{equation} 
between the space $C^{\infty}(P^n\mathbb{H})$ and the
space $\Gamma_{\mathcal{F}}(L\otimes\sqrt{K_{\Ha}^{\mathcal{F}}})$.  
Again here note that $\nabla_X(\bs{s}_{\Ha}) = 0 $ 
for $X\in \Gamma (\mathcal{F})$.
So we can define the inner
product on the space 
$\Gamma_{\mathcal{F}}(L\otimes\sqrt{K_{\Ha}^{\mathcal{F}}})$ by
\begin{equation}
\int_{P^n\mathbb{H}}g_1(P)\overline{g_2(P)}\bs{v}_{\Ha}
\end{equation}
through the above identification \eqref{eq:3.11}, since we assumed
that $\inner{\bs{s}_{\Ha}}{\bs{s}_{\Ha}}_L\equiv 1$ and 
$\sqrt{\bs{v}_{\Ha}}\otimes \sqrt{\bs{v}_{\Ha}} = \bs{v}_{\Ha}$.

By \eqref{eq:bH} we define the pairing 
$< \sqrt{\pi_{\Ha}\,^*(\bs{v}_{\Ha})}, \sqrt{\sigma_{\Ha}} >$
as
\begin{equation}
<\sqrt{\pi_{\Ha}\,^*(\bs{v}_{\Ha})}, \sqrt{\sigma_{\Ha}}> = 
\sqrt{|\bs{b}_{\Ha}|}\,\norm{A}^{\frac{1}{2}},
\end{equation}       
then finally 
we define the pairing between the spaces 
$\Gamma_{\mathcal{F}}(L\otimes\sqrt{K_{\Ha}^{\mathcal{F}}})$ 
and $\Gamma_{\mathcal{G}}(L\otimes\sqrt{K_{\Ha}^{\mathcal{G}}})$ as follows:
Let $\varphi=
\pi_{\Ha}\,^*(f)\cdot\bs{s}_{\Ha}\otimes\sqrt{\pi_{\Ha}\,^*(\bs{v}_{\Ha})}$ and
$\psi=g \cdot \bs{t}_{\Ha}\otimes\sqrt{\sigma_{\Ha}}$, 
then we define
\begin{equation}
\pair{\varphi}{\psi}\equiv\int_{\widetilde{\mathbf{E}_{\Ha}}}
\pi_{\Ha}\,^*(f)(A)\,g(A)\,(\bs{t}_{\Ha},\bs{s}_{\Ha})_L\,
<\sqrt{\pi_{\Ha}\,^*(\bs{v}_{\Ha})}, \sqrt{\sigma_{\Ha}}>\, \Omega_{\Ha}
\end{equation}
\begin{equation*}
=\int_{\widetilde{\mathbf{E}_{\Ha}}}
\pi_{\Ha}\,^*(f)g(A)e^{-2^{\frac 14}\pi
\sqrt{\norm{A}}}\sqrt{|\bs{b}_{\Ha}|}\norm{A}^{\frac{1}{2}}\,\Omega_{\Ha}.
\end{equation*}

The pairing above defines
an operator $T : g \mapsto T(g) \in C^{\infty}(P^n\Ha)$
for $g\in C^{\infty}(\widetilde{\mathbf{E}_{\Ha}})$ satisfying
a suitable integrability condition, that is we have 
\begin{equation*}
\int_{P^n\Ha}f\,T(g)\bs{v}_{\Ha}
=\pair{\pi^*_{\Ha}(f)\cdot\bs{s}_{\Ha}\otimes\sqrt{\pi_{\Ha}\,^*(\bs{v}_{\Ha})}}
{\,g \cdot \bs{t}_{\Ha}\otimes\sqrt{\sigma_{\Ha}}}. 
\end{equation*}
In fact we have
\begin{equation}
T(g)\bs{v}_{\Ha}=(\pi_{\Ha})_{*}
(ge^{-2^{\frac 14}\pi\sqrt{\norm{A}}}\sqrt{|\bs{b}_{\Ha}|}
\norm{A}^{\frac{1}{2}}\Omega_{\Ha}),
\end{equation}
the fiber integration of the map $\pi_{\Ha}$.

Let $\Delta^{\Ha}$ be the Laplacian on $P^n\Ha$ with respect to the
Riemannian metric defined in $\S\,2$
and let
$\{e^{-\sqrt{-1}t\sqrt{\Delta^{\Ha}+(2n+1)^2}}\}_{t\in\mathbb{R}}$ 
be the one parameter family of unitary and elliptic Fourier integral operators 
generated by the square root
of the operator $\Delta^{\Ha}+ (2n+1)^2$.  

The bicharacteristic flow $\{\sigma_t\}_{t\in\R}$ of the operator 
$\sqrt{\Delta^{\Ha}+(2n+1)^2}$  is the Hamilton flow whose 
Hamiltonian is the square root of the metric function and can be expressed
as follows:
\begin{prop}
\begin{equation}
\sigma_t : \widetilde{\mathbf{E}_{\Ha}} \rightarrow
\widetilde{\mathbf{E}_{\Ha}},\,\,
\sigma_t(A) = e^{-2\sqrt{-1}t}\,A.
\end{equation}
\end{prop}

\begin{proof}
Let $X_g$ be the vector field on $\widetilde{\mathbf{E}_{\Ha}}$ defined
by the flow $\sigma_t(A)= e^{-2\sqrt{-1}t}\,A$, and let 
$h(A) = 2^{-3/4}\sqrt{\norm{A}}$(= square root of the metric function
on $TP^n\Ha$). Then by Proposition \ref{sec1:kaehler} 
we have for any vector field
$Y$ on $\widetilde{\mathbf{E}_{\Ha}}$ 
\begin{equation*}
\omega_{\Ha}(Y,X_g)= 
2^{1/4}\sqrt{-1}\,\overline{\partial}{\partial}\sqrt{\norm{A}}(Y,X_g)= Y(h),
\end{equation*}
which proves the proposition.
\end{proof}

The geodesic flow restricted to the unit tangent sphere bundle of
$P^n\Ha$ coincides with the Hamilton flow above under the
identification of the tangent bundle and the cotangent bundle by
the Riemannian metric.

{}From the theory of Fourier integral operators, we know that 
for each $t\in\R$ the graph of the symplectic isomorphism $\sigma_t$
is the canonical relation of the Fourier integral operator
$e^{-\sqrt{-1}t\sqrt{\Delta^{\Ha}+(2n+1)^2}}$. 
This correspondence can be seen as a quasi-classical approximation. 
An opposite is interpreted as a kind of quantization. Here we have a
kind of exact
quantization of the flow $\{\sigma_t\}$ by making use of the operator $T$:  
since $\sigma_t^{\,*}(\sigma_{\Ha}) = e^{-2\sqrt{-1}t(2n+1)}\sigma_{\Ha}$ 
and $\sigma_t^{\,*}(\bs{t}_{\Ha}) = \bs{t}_{\Ha}$, we should regard the
action of the flow $\{\sigma_t\}$  on the Hilbert space 
$\mathfrak{h}^{\mathcal{G}}_{\Ha}$ such that
$(\sigma_t)^*(f\cdot\bs{t}_{\Ha}\otimes\sqrt{\sigma_{\Ha}})$ =
$e^{-\sqrt{-1}t(2l+2n+1)}\,f\cdot\bs{t}_{\Ha}\otimes\sqrt{\sigma_{\Ha}}$ 
for each $f\in\mathcal{P}^{\Ha}_l$.
Then we have  
\begin{thm}\label{th:3.2}
\begin{enumerate}
\item The operator $T$ is an isomorphism between the spaces
$\mathfrak{h}_{\Ha}^{\mathcal G}$ and $L_2(P^n\Ha)$.
\item The following diagram is commutative
\begin{equation}
\begin{CD}
\mathfrak{h}_{\Ha}^{\mathcal G} @>{\sigma_t^*}>>
\mathfrak{h}_{\Ha}^{\mathcal G}\\
@V{T}VV  @VV{T}V\\
L_2(P^n\Ha) @>>{e^{-t\sqrt{-1}\sqrt{\Delta^{\Ha}+(2n+1)^2}}}>
L_2(P^n\Ha).
\end{CD}
\end{equation}
\end{enumerate}
\end{thm}

\begin{rem}
The Hamilton flow $\{\sigma_t\}$ is periodic with the period $\pi$.
However, the action $\{\sigma_t^*\}$ on $\mathfrak{h}^{\mathcal{G}}_{\Ha}$
is periodic with the period $2\pi$ and the one parameter family of
unitary operators $\{e^{-t\sqrt{-1}\sqrt{\Delta^{\Ha}+(2n+1)^2}}\}$
also has the period $2\pi$.
\end{rem} 

We will prove this theorem by explicitly calculating the
norms of the operator $T$ on each $\mathcal{P}_l^{\Ha}$
( Proposition \ref{prop:3.8} ).

Let us denote by $\mathcal{S}_k\,(k=0,1,\dots)$ the space of
harmonic polynomials on $\Ha^{n+1}\cong\R^{4n+4}$ of
degree $k$, and by $\mathcal{S}_k^0$ the subspace of
$\mathcal{S}_k$ consisting of those polynomials which are
invariant under the action of $\mathit{Sp}(1)$ from the right.
The $l$-th eigenspace $H_l$ of the Laplacian $\Delta^{\Ha}$ on
$P^n\Ha$ with the eigenvalue $\lambda_l=4l(2n+1+l)$ is isomorphic
to $\mathcal{S}_{2l}^0$ by $\pi^* : H_l\cong\mathcal{S}_{2l}^0$($\pi :
S^{4n+3}$ $\rightarrow$ $P^n{\Ha}$),
and $\mathcal{S}_{2l+1}^0=\{0\}$.  It is known that the dimension of
$H_l$ is given by
\begin{equation}
\dim \,H_l = \frac{2n}{2n+1}\frac{l+1}{l+2n}\cdot (2l+2n+1)\cdot
\left(\frac{\Gamma(l+2n+1)}{\Gamma(2n+1)\Gamma(l+2)}\right)^2,
\end{equation}
and 
\begin{equation}
\overline{\sum_{l=0}^{\infty} H_l} = L_2(P^n{\Ha}).
\end{equation}

Here we note the important fact that
for each fixed $z \in \Ha^{n+1}\otimes\C$ satisfying $\Einner{z}{z}=0$,
the polynomial $\Einner{p}{z}^l$ of $p\in\Ha^{n+1}$ $\cong$ $\R^{4n+4}$
is a harmonic polynomial, and likewise we can prove
\begin{prop}\label{prop:3.4}
For each fixed $A\in \widetilde{\mathbf{E}_{\Ha}}$, the polynomial 
$\Einner{\pi(p)}{A}_{\C}^{\,\,l}$ =
$\Einner{(\theta(p_i)p_j)}{A}_{\C}^{\,\,l}$ ($l=0,1,\dots$) 
on $\Ha^{n+1}\cong \R^{4n+4}$ is a harmonic polynomial and $Sp(1)$-invariant.
\end{prop}

\begin{proof}
It is enough to prove that the polynomial 
$\Einner{\pi(p)}{A}_{\C}^{\,\,l}$
is harmonic for a particular $A \in \widetilde{\mathbf{E}_{\Ha}}$
because of the two point homogeneity of $P^n\Ha$ by the isometry group
$Sp(n+1)$.  Hence if we take a point $A\in\widetilde{\mathbf{E}_{\Ha}}$, 
\begin{equation*}
A = \rho
\left(
\left(\begin{array}{ccccc}
1 & \sqrt{-1} & 0 &\cdots & 0\\
\sqrt{-1} & -1 &0 &\cdots & 0\\
0 & 0 &&&\\
\vdots & \vdots &&\bigzerol&\\
0 & 0 &&&
\end{array}\right)
\right),
\end{equation*}
then we have the result by an explicit calculation. 
\end{proof}

Let $\mathcal{A}_l$ be a map:
\begin{equation}
\mathcal{A}_l:\; H_l\rightarrow \mathcal{P}_l^{\Ha}
\end{equation}
given by
\begin{equation}
\mathcal{A}_l(\varphi)(A)=\int_{P^n\Ha}\varphi(P)
\Einner{P}{A}_{\C}^l\bs{v}_{\Ha}.
\end{equation}
For $f\in\mathcal{P}^{\Ha}_l$ let 
$\mathcal{B}_l(f)\in C^{\infty}(P^n\Ha)$ 
be given by the integral
\begin{equation}
\mathcal{B}_l(f)(P)=\int_{\widetilde{\mathbf{E}_{\Ha}}}f(A)
\overline{\Einner{P}{A}_{\C}^l} e^{-2\sqrt[4]{2}\pi\sqrt{\norm{A}}}
\sqrt{|\bs{a}_{\Ha}|}\norm{A}^{n+1}\Omega_{\Ha}.
\end{equation}
Then Proposition \ref{prop:3.4} guarantees that 
$\mathcal{B}_l$ maps $\mathcal{P}^{\Ha}_l$ into $H_l$ and
we have 
\begin{enumerate}
\item $\mathcal{B}_l$ commutes with the action of $Sp(n+1)$ on $\mathcal{P}_l$
and $H_l$.
\item $\displaystyle\int_{P^n\Ha}(\mathcal{B}_l\circ A_l)(\varphi)(P)\cdot
\overline{(\varphi)(P)}\bs{v}_{\Ha}\\
=\int |A_l(\varphi)(P)|^2\cdot
e^{-2\sqrt[4]{2}\cdot\pi\cdot\sqrt{\norm{A}}}\sqrt{|\bs{a}_{\Ha}|}
\norm{A}^{n+1}\Omega_{\Ha}$.
\end{enumerate}
Of course the operator $\mathcal{A}_l$ also commutes with the action
of $Sp(n+1)$ on $H_l$ and $\mathcal{P}^{\Ha}_l$, so that
the operator 
$\mathcal{B}_l\circ\mathcal {A}_l$ 
is a positive constant multiple (= $\bs{b}_l$) of the identity operator.
This constant $\bs{b}_l$ satisfies
\begin{multline}
\int_{\widetilde{\mathbf{E}_{\Ha}}}
|\mathcal{A}_l(\varphi)(A)|^2\;e^{-2\sqrt[4]{2}
\pi\sqrt{\norm{A}}}\sqrt{|\bs{a}_{\Ha}|}
\norm{A}^{n+1}\Omega_{\Ha}\\
=\bs{b}_l\int_{P^n\Ha}|\varphi|^2\;\bs{v}_{\Ha}.
\qquad\qquad\qquad\qquad\qquad
\end{multline}
To determine $\bs{b}_l$ we put
\begin{equation}
G_l(P,P')=\int_{\widetilde{\mathbf{E}_{\Ha}}}\Einner{P}{A}_{\C}^l
\overline{\Einner{P'}{A}_{\C}^l}e^{-2\sqrt[4]{2}
\pi\sqrt{\norm{A}}}\sqrt{|\bs{a}_{\Ha}|}
\norm{A}^{n+1}\Omega_{\Ha},
\end{equation}
then
\begin{equation}
\int_{P^n\Ha}G_l(P,P')\varphi(P)\bs{v}_{\Ha}=\bs{b}_l\varphi(P'),\quad
\varphi\in H_l.
\end{equation}
Now from the invariance $G_l(gP{\Trans{g}},g{P'}{\Trans{g}})=G_l(P,P')$
for any $g\in Sp(n+1)$, we have
\begin{equation}
\int_{P^n\Ha}G(P,P)\bs{v}_{\Ha}=\bs{b}_l\dim H_l=G(P_0,P_0)\vol(P^n\Ha)
\end{equation}
for any fixed $P_0\in P^n\Ha$, and we have
\begin{multline}
\bs{b}_l\dim H_l\\
=\int_{\widetilde{\mathbf{E}_{\Ha}}}\int_{P^n\Ha}\left|\Einner{P}{A}_{\C}
\right|^{2l}e^{-2\sqrt[4]{2}\pi\sqrt{\norm{A}}}\sqrt{|\bs{a}_{\Ha}|}
\norm{A}^{n+1}\bs{v}_{\Ha}\wedge\Omega_{\Ha}.
\end{multline}
Now the integral\begin{sloppypar}
\begin{equation}
\norm{A}^{-2l}\int_{P^n\Ha}\left|\Einner{P}{A}_{\C}\right|^{2l}\bs{v}_{\Ha}
=I_l
\end{equation}
is independent of $A\in\widetilde{\mathbf{E}_{\Ha}}$, so that we have
\begin{eqnarray*}
\bs{b}_l\dim H_l
&=&\int_{\widetilde{\mathbf{E}_{\Ha}}}I_l\norm{A}^{n+1+2l}
\sqrt{|\bs{a}_{\Ha}|}e^{-2\sqrt[4]{2}\pi\sqrt{\norm{A}}}\Omega_{\Ha}\\
&=& I_l\sqrt{|\bs{a}_{\Ha}|}2^{\frac 32 (n+1+2l)}\int_{T_0^*P^n\Ha}
\norm{q}^{2n+2+4l}e^{-4\pi\norm{q}}\Omega_{\Ha}\\
&=& I_l\sqrt{|\bs{a}_{\Ha}|}2^{\frac 32 (n+1+2l)}\\
&\qquad&\times\int_{S^{4n-1}}\int_0^{\infty}
r^{2n+2+4l}e^{-4\pi r}r^{4n-1}\,dr\wedge\bs{v}_{S^{4n-1}}
{\cdot}\vol(P^n\Ha)\\
&=& I_l\sqrt{|\bs{a}_{\Ha}|}2^{\frac 32 (n+1+2l)}\frac{1}{(4\pi)^{6n+4l+2}}\\
&\qquad&\times\Gamma(6n+4l+2)\vol(S^{4n-1})\vol(P^n\Ha).
\end{eqnarray*}
\end{sloppypar}
                  
Next we calculate the constant $I_l$. To this purpose we choose a point
$(P_0,Q_0)$ from $\mathbf{E}_{\Ha}$:
\begin{equation*}
P_0=\left(\begin{array}{cccc}
1 & 0 & \cdots & 0\\
0 &&&\\
\vdots &&\bigzerol&\\
0 &&&
\end{array}\right)
\;\,\,\text {and}\; \,\,
Q_0=\left(\begin{array}{ccccc}
0 & 1 & 0 &\cdots & 0\\
1 & 0 &0 &\cdots & 0\\
0 & 0 &&&\\
\vdots & \vdots &&\bigzerol&\\
0 & 0 &&&
\end{array}\right). 
\end{equation*}
When we put
$A_0=\tau_{\Ha}(P_0,Q_0)$,
then
\begin{equation}
\left|\Einner{P}{A_0}_{\C}\right|^2=(\norm{p_0}^2-\norm{p_1}^2)^2
+4\left|\Einner{p_0}{p_1}\right|^2,
\end{equation}
where $P=(p_i\theta(p_j))$, $(p_0,\dots,p_n)\in S^{4n+3}\subset
\Ha^{n+1}$. So
\begin{multline}\label{eq:3.30}
I_l=\norm{A_0}^{-2l}\int_{P^n\Ha}\left|\Einner{P}{A_0}_{\C}\right|^{2l}
\bs{v}_{\Ha}\\
= \norm{A_0}^{-2l} \int_{P^n\Ha} \left|\Einner{P}{A_0}_{\C}\right|^{2l}
\cdot \frac{1}{2\pi^2}\pi_*(\bs{v}_S)\\
=\frac{(2\sqrt{2})^{-2l}}{2\pi^2}\int_{S^{4n+3}}
((\norm{p_0}^2-\norm{p_1}^2)^2+4\left|\Einner{p_0}{p_1}\right|^2)^l
\bs{v}_{S}.
\end{multline}

For the determination of the last integral on the sphere $S^{4n+3}$,
we define a coordinate transformation $\Phi$:
\begin{equation*}
\Phi:\;\Ha^2\times D_1
\rightarrow\Ha^{n+1},(y,\bs{x})\mapsto(\sqrt{1-|\bs{x}|^2}y,\bs{x}),
\end{equation*}
where $D_1=\{\bs{x}\in\Ha^{n-1}|\;\norm{\bs{x}}<1\}$.
Then we can separate the variables into two parts $(y,x)$ in the last 
expression of \eqref{eq:3.30} so that the
constant $I_l$ is equal to the following integral:
\begin{multline}\label{eq:3.31}
I_l = \frac{(2\sqrt{2})^{-2l}}{2\pi^2}\int_{\norm{\bs{x}}<1}(1-|\bs{x}|^2)^{3+2l}
\,d\bs{x}\\
\times\int_{S^7}((\norm{y_0}^2-\norm{y_1}^2)^2 + 4
\left|\Einner{y_0}{y_1}\right|^2)^l \bs{v}_{S^7},
\end{multline}
where $d\bs{x}$ is the Lebesgue measure on $\Ha^{n-1}\cong\R^{4n-4}$
and $\bs{v}_{S^7}$ is the volume element on the unit sphere 
$S^7\subset\Ha^2\cong\R^8$.

Let $\pi:\;\Ha^2\cong\R^8\rightarrow\R\times\Ha\cong\R^5$ be the map
\begin{center}
\quad $\pi(y_0,y_1)=(\norm{y_0}^2,y_0\theta(y_1))=(t_0,t_1,t_2,t_3,t_4)$.
\end{center}
Then $\pi$ realizes the Hopf fibration\\
\begin{equation*}
\pi:\;S^7\rightarrow S^4=\left\{(t_0,\dots,t_4)\left|\;
\left(t_0-\frac 12\right)^2+\sum_{i=1}^4 t_i^2=\frac 14\right.\right\},
\end{equation*}
and we have a reduction of the integral on $S^7$ 
in the above formula \eqref{eq:3.31} 
by the fiber integration to an integral on $S^4$:
\begin{eqnarray*}
&&\int_{S^7}\left((\norm{y_0}^2-\norm{y_1}^2)^2+4\left|\Einner{y_0}{y_1}
\right|^2\right)^l\bs{v}_{S^7}\\
&=& \int_{S^4}\pi_{*}\left(\left((\norm{y_0}^2-\norm{y_1}^2)^2+
4\left|\Einner{y_0}{y_1}\right|^2\right)^l\bs{v}_{S^7}\right)\\
&=&\frac{\vol(S^3)}{2^4}\int_{S^4}(t_0^2+t_1^2)^l\bs{v}_{S^4}\\
&=&\frac{\vol(S^3)}{2^4}2\pi\int_{t_2^2+t_3^2+t_4^2<1}(1-t_2^2-t_3^2-t_4^2)^l\,
dt_2dt_3dt_4\\
&=&\frac{\vol(S^3)\vol(S^2)}{2^4}\cdot {2\pi}\cdot
\frac{\Gamma(l+1)\Gamma\left(\frac 32\right)}{2\Gamma\left(l+1+\frac 32\right)}.
\end{eqnarray*}
Hence we have 
\begin{multline}
I_l = \frac{(2\sqrt{2})^{-2l}}{2\pi^2}
\cdot\vol(S^{4n-5})
\frac{\Gamma(2l+4)\Gamma(2n-2)}{2\Gamma(2l+2n+2)}\\
\times\frac{\vol(S^3)\vol(S^2)}{2^4}\cdot 2\pi
\cdot\frac{\Gamma(l+1)\Gamma\left(\frac 32\right)}
{2\Gamma\left(l+1+\frac 32\right)}.
\end{multline}

Finally the constant $\bs{b}_l$ is expressed as
\begin{prop}\label{prop:3.5}
\begin{align}
\bs{b}_l&=\frac{\sqrt{|\bs{a}_{\Ha}|}}{\dim H_l}2^{\frac 34 (2n+2+4l)}
\frac{\vol(S^{4n-1})\vol(P^n\Ha)}{(4\pi)^{6n+4l+2}}\Gamma(6n+4l+2)\\
&\qquad\qquad\times\frac{(2\sqrt{2})^{-2l}}{2\pi^2}\frac{\vol(S^3)\vol(S^2)}{2^4}
{2\pi}\frac{\Gamma(l+1)\Gamma\left(\frac 32\right)}
{2\Gamma\left(l+1+\frac 32\right)}\notag\\
&\qquad\qquad\times\frac{\vol(S^{4n-5})\Gamma(2l+4)\Gamma(2n-2)}
{2\Gamma(2l+2n+2)}\notag\\
&=\sqrt{|\bs{a}_{\Ha}|}\cdot 2^{-n/2+1}\cdot \pi^{-4l-3}\cdot
\frac{1}{(2l+2n+1)^2}\notag\\
&\,\,\,\times\frac{\Gamma(l+1)^2\Gamma(l+2)^2}
{\Gamma(l+n+1/2)\Gamma(l+n+1)\Gamma(l+2n)\Gamma(l+2n+1)}\notag\\
&\,\,\,\,\times\Gamma(l+\frac{6n+2}{4})\Gamma(l+\frac{6n+3}{4})
\Gamma(l+\frac{6n+4}{4})\Gamma(l+\frac{6n+5}{4})\notag
\end{align}
\end{prop}
As a result we proved 
\begin{prop}
The operator $\frac{1}{\sqrt{\bs{b}_l}}\mathcal{A}_l$ is
a unitary isomorphism between the spaces $H_l$ and
$\mathcal{P}_l^{\Ha}$, and $\frac{1}{\sqrt{\bs{b}_l}}\mathcal{B}_l$ is
the inverse.
\end{prop}

Next we consider the operator $T\circ\mathcal{A}_l$ $:$
$H_l \rightarrow C^{\infty}(P^n\Ha)$.  

Let $\varphi\in H_l$, then
\begin{eqnarray*}
&& (T\circ\mathcal{A}_l(\varphi))\bs{v}_{\Ha}\\
&=& (\pi_{\Ha})_* \left(\left(\int_{P^n\Ha}\varphi(P)
\Einner{P}{A}_{\C}^l
\bs{v}_{\Ha}\right)e^{-\sqrt[4]{2}\cdot\pi\cdot\sqrt{\norm{A}}}
\sqrt{|\bs{b}_{\Ha}|}\norm{A}^{\frac{1}{2}}\Omega_{\Ha}\right),\\
\end{eqnarray*}
so put $(\pi_{\Ha})_*\left(\Einner{P}{A}_{\C}^l
e^{-\sqrt[4]{2}\cdot\pi\cdot\sqrt{\norm{A}}}
\sqrt{|\bs{b}_{\Ha}|}\norm{A}^{\frac{1}{2}}\Omega_{\Ha}\right)
=K_l(P,P')\bs{v}_{\Ha}$, where $\pi_{\Ha}(A)=P'$ (exactly
$P'=\pi_{\Ha}\circ\tau_{\Ha}^{-1}(A)$). 
Then 
\begin{equation*}
T\circ\mathcal{A}_l(\varphi)(P') 
= \int_{P^n\Ha}K_l(P,P')\varphi(P)\bs{v}_{\Ha}.
\end{equation*}
Since $K_l(gP{\Trans{g}},\,g{P'}{\Trans{g}})=K_l(P,\,P')$ 
for any $g\in\mathit{Sp}(n+1)$ and so $K_l(P,\,P')$ = $K_l(P',\,P)$
(since $P^n\Ha$ is a symmetric space),
we know that for each fixed $P'$ (resp. $P$) the polynomial 
$K_l(\pi(p)\,,\,P')$ (resp. $K_l(P,\,\pi(p')$) on $\Ha^{n+1}$ 
is a harmonic polynomial. Hence the operator 
$T\circ\mathcal{A}_l$ maps $H_l$ into itself, so that 
$T\circ\mathcal{A}_l:\;H_l\rightarrow H_l$ is an
intertwining operator of the irreducible unitary
representation of $\mathit{Sp}(n+1)$ on $H_l$. Hence 
we can put
\begin{equation*}
T\circ\mathcal{A}_l =\bs{a}_l Id,
\end{equation*}
on $H_l$ with a suitable constant $\bs{a}_l$.

\begin{prop}
\begin{multline}
\bs{a}_l=\frac{(2\sqrt{2})^{\frac{1}{2}}}{\dim H_l}
\frac{\Gamma(2l+4n+1)}{(2\pi)^{2l+4n+1}}
\mathit{vol}(S^{4n-1})\mathit{vol}(P^n\Ha)\cdot\sqrt{|\bs{b}_{\Ha}|}\\
=\sqrt{|\bs{b}_{\Ha}|}\cdot\frac{2^{3/4}}{\pi^{2l+3/2}}\cdot
\frac{\Gamma(l+1)\Gamma(l+2)\Gamma(l+2n+1/2)}{(2l+2n+1)\,\Gamma(l+2n)}
\end{multline}
\end{prop}

\begin{proof}
Since 
\begin{equation*}
\int_{P^n\Ha}K_l(P,P')\varphi(P)\bs{v}_{\Ha}=\bs{a}_l\varphi(P'),
\end{equation*}
and by the invariance of the kernel $K_l(P,P')$ $=$
$K_l(gP{\Trans{g}},g{P'}{\Trans{g}})$,
we have
\begin{equation*}
\int_{P^n\Ha}K_l(P,P)\bs{v}_{\Ha}=\bs{a}_l\dim H_l
=K_l(P_0,P_0)\mathit{vol}(P^n\Ha),\quad P_0\in P^n\Ha.
\end{equation*}

Let $\tau_{\Ha}(P_0,Q)=A$, then
\begin{eqnarray*}
\Einner{P_0}{A}_{\C}
&=&\Einner{P_0}{(\norm{Q}^2P_0-Q^2)\otimes 1+\dfrac{\norm{Q}}{\sqrt{2}}
\otimes\sqrt{-1}}_{\C}\\
&=&\frac 12 \norm{Q}^2.
\end{eqnarray*}
Hence
\begin{multline*}
K_l(P_0,P_0)\bs{v}_{\Ha}\\
=(\pi_{\Ha})_*\left(\left(\frac 12 \norm{Q}^2
\right)^le^{-\sqrt{2}\cdot\pi\cdot\norm{Q}}
\sqrt{|\bs{b}_{\Ha}|}({\sqrt{2}\norm{Q}^2})^{\frac{1}{2}}\Omega_{\Ha}
\right).
\end{multline*}
Since the integrand on the fiber $\pi_{\Ha}^{-1}(P_0)$ is only
a function of the norm $\norm{Q}$, we have
\begin{align*}
K_l(P_0,P_0)
&=\int_0^{\infty}t^{2l+1}e^{-2\pi t}t^{4n-1}\,dt\;\mathit{vol}
(S^{4n-1})\sqrt{|\bs{b}_{\Ha}|}(2\sqrt{2})^{\frac{1}{2}}\\
&=\frac{(2\sqrt{2})^{\frac{1}{2}}}{(2\pi)^{2l+4n+1}}
\Gamma(2l+4n+1)\vol(S^{4n-1})\sqrt{|\bs{b}_{\Ha}|}.
\end{align*}
Note here that $\norm{Q}^2=2\norm{q}^2$ for $\inner{p}{q}_{\Ha}=0,
\alpha(p,q)=(P,Q)$ and $2\norm{Q}^4=\norm{A}^2$, $A=\tau_{\Ha}
(P,Q)$.
\end{proof}

The proof of the commutativity in Theorem \ref{th:3.2} is included in 
the above arguments and now we have proved Theorem \ref{th:3.2}. 
We restate the results more precisely in  
\begin{prop}\label{prop:3.8}
\begin{enumerate}
\item The operator $T$ is a constant multiple of a unitary
operator on each $\mathcal{P}_l^{\Ha}$.
\item The norm of $T$ on $\mathcal{P}_l^{\Ha}$ is given by
\begin{align}
\norm{T_{|\mathcal{P}^{\Ha}_{l}}}&=\frac{\bs{a}_l}{\sqrt{\bs{b}_l}}\\
&=\frac{\sqrt{\bs{b}_{\Ha}}}{\sqrt[4]{\bs{a}_{\Ha}}}\cdot
2^{(n+1)/4}\cdot\frac{\Gamma(l+2n+1/2)}{\Gamma(l+2n)}\notag\\
&\times\sqrt{\frac{\Gamma(l+n+1/2)\Gamma(l+n+1)\Gamma(l+2n)\Gamma(l+2n+1)}
{\Gamma(l+\frac{6n+2}{4})\Gamma(l+\frac{6n+3}{4})\Gamma(l+\frac{6n+4}{4})
\Gamma(l+\frac{6n+5}{4})}}\notag
\end{align}
\item By an asymptotic property of products of the Gamma function (see \cite{Ra2}) 
\begin{equation}
\lim_{l\rightarrow\infty} \norm{T_{|\mathcal{P}^{\Ha}_l}} 
=\frac{\sqrt{\bs{b}_{\Ha}}}{\sqrt[4]{\bs{a}_{\Ha}}}2^{\frac{n+1}{4}}
= \frac{\sqrt{2}}{\pi},
\end{equation}
which also proves that the operator $T$ is an isomorphism between
the two Hilbert spaces $\mathfrak{h}^{\mathcal{G}}_{\Ha}$ and
$L_2(P^n\Ha)$.
\end{enumerate}
\end{prop}


\section{Quantization operator II}

In this section we describe 
another quantization operator $\tilde{T}:
\mathfrak{h}^{\mathcal{G}}_{\Ha} \rightarrow L^2(P^n\Ha)$. 
We construct it by making use of Proposition \ref{prop:1.7}
(see ~\cite{FY} for the case of the complex projective space).

Based on the data explained in $\S\, 2$ we can construct a
quantization operator for the case of the sphere $S^{4n+3}$(for
details see ~\cite{Ra2}). We denote this operator by $T^S$. 
It is expressed as a fiber integration:
\begin{equation}
T^S(g)(\pi_S(B))\bs{v}_S = 
(\pi_S)_*(g(B)\cdot e^{-\pi\norm{B}}\sqrt{|\bs{b}_S|}\norm{B}^{-1/2}\Omega_S),
\end{equation}
where $g$ is a function on $T^*_0S^{4n+3} \cong \widetilde{\mathbf{E}_S}$ 
satisfying a suitable integrability condition.
Then the operator 
 $\tilde{T}$ $:$ $\sum_l\mathcal{P}^{\Ha}_l$  
$\rightarrow$ $C^{\infty}(P^n\Ha)$ 
is defined as follows:
let $f\in \mathcal{P}^{\Ha}_l$, 
then $T^S(\beta^*(f))$ is $Sp(1)$ invariant, so that it can descend
to $P^n\Ha$ and we denote the resulting function by
$\tilde{T}(f)(P)$, or
\begin{align}
&\tilde{T}(f)\bs{v}_{\Ha}\\ 
&=\frac{1}{2\pi^2}(\pi)_*(T^S(\beta^*(f))\bs{v}_S)\notag
\end{align}
Note that $\beta^*(f)$ is a polynomial of degree $2l$ on  
$\underbrace{M(2,\C)\times\dots\times M(2,\C)}_{n+1}$. 
By the same arguments for the operator $T\circ \mathcal{A}_l$, we see
$\tilde{T}\circ\mathcal{A}_l$ $=$ $\bs{c}_l\cdot Id$ 
on each subspace $H_l$ with a constant $\bs{c}_l$.

\begin{prop}
\begin{align*}
\quad\bs{c}_l\\
&=\sqrt{|\bs{b}_S|}\frac{\vol(P^n\Ha)}{\dim H_l}\frac{\sqrt{\pi}}{4}
\frac{\vol(S^2)\vol(S^{4n-1})}{l+2n+1/2}\frac{\Gamma(l+2n)}{\Gamma(l+2n+1/2)}\\
&\qquad\qquad\times\frac{1}{\sqrt{2}}\frac{\Gamma(2l+4n+5/2)}
{(2\pi)^{2l+4n+5/2}}\\
&=\frac{\sqrt{|\bs{b}_S|}}{2\sqrt{2}\pi^{\frac{3}{2}}}\cdot
\frac{1}{\pi^{2l}}\cdot\frac{(l+2n+1/4)(l+2n+3/4)}{(l+2n+1/2)^2}
\cdot\Gamma(l+1)\Gamma(l+2)\\
&\qquad\qquad\times\frac{\Gamma(l+2n+1/4)\Gamma(l+2n+3/4)}
{\Gamma(l+2n+1)\Gamma(l+2n+1/2)}
\end{align*}
\end{prop}

\begin{proof}
Let $\varphi \in H_l$, then 
\begin{align*}
&(\tilde{T}^S\circ
\mathcal{A}_l)(\varphi)(P')\bs{v}_{\Ha}\\
& =\frac{1}{2\pi^2}
(\pi\circ\pi_S)_*\left(\left(\int_{P^n\Ha}\varphi(P)
\Einner{P}{\beta (B)}_{\C}^{\,l}\bs{v}_{\Ha}\right)\cdot
e^{-\pi\norm{B}}\sqrt{|\bs{b}_S|}\norm{B}^{-1/2}\Omega_S\right)\\
&=\int_{P^n\Ha}\varphi(P)L_l(P,P')\bs{v}_{\Ha}=\bs{c}_l\varphi(P),
\end{align*}
where 
$L_l(P,P') = \tilde{T}(\Einner{P}{\beta(B)}_{\C}^{\,\,l})$.
Of course, the fiber integration is taken with respect to the variable $B$,
and we put
$P'= \pi\circ\pi_S(B)$.
As before, to determine the constant $\bs{c}_l$ it is enough to
calculate the integral
\begin{equation*}
\int_{P^n\Ha}L_l(P,P)\bs{v}_{\Ha} = \bs{c}_l \dim H_l. 
\end{equation*}
Because of the invariance of $L_l$, we have for any point
$P\in P^n\Ha$
\begin{align*}
&\bs{c}_l\frac{\dim H_l}{\vol(P^n\Ha)}\,\bs{v}_{\Ha}
= L_l(P,P)\bs{v}_{\Ha}\\
& =\frac{1}{2\pi^2} (\pi\circ\pi_S)_*(\Einner{P}{\beta(B)}_{\C}^{\,l}
e^{-\pi\norm{B}}\sqrt{|\bs{b}_S|}\norm{B}^{-1/2}\Omega_S)\\
&= \frac{1}{2\pi^2}(\pi\circ\pi_S)_*(((1/2)\norm{Q}^2)^le^{-\pi\norm{B}}
\sqrt{|\bs{b}_S|}\norm{B}^{-1/2}\Omega_S)\\
&= \frac{1}{2\pi^2}\sqrt{|\bs{b}_S|}
(\pi\circ\pi_S)_*((\norm{q}^2-\norm{(q,p)_{\Ha}}^2)^l
e^{-2\pi\norm{q}}(2\norm{q})^{-1/2}\Omega_S),
\end{align*}
where we put $B = \norm{q}\rho(p) + \rho(q)\sqrt{-1}$ and
$\pi(p)=P$.  Then we have
\begin{align*}
&L_l(P,P)\\
&\equiv \sqrt{|\bs{b}_S|}\cdot\int_{\R^{4n+3}}
(\norm{x}^2 - (x_1^2+x_2^2+x_3^2))^l\cdot 
e^{-2\pi\norm{x}}(2\norm{x})^{-1/2}dx\\
&=\sqrt{|\bs{b}_S|}\frac{\sqrt{\pi}}{4}
\frac{\vol(S^2)\vol(S^{4n-1})}{l+2n+1/2}
\frac{\Gamma(l+2n)}{\Gamma(l+2n+1/2)}
\frac{1}{\sqrt{2}}\frac{\Gamma(2l+4n+5/2)}{(2\pi)^{2l+4n+5/2}}.
\end{align*}
Finally we have
\begin{align*}
&\quad\bs{c}_l\\
&=\sqrt{|\bs{b}_S|}\frac{\vol(P^n\Ha)}{\dim H_l}\frac{\sqrt{\pi}}{4}
\frac{\vol(S^2)\vol(S^{4n-1})}{l+2n+1/2}\frac{\Gamma(l+2n)}{\Gamma(l+2n+1/2)}\\
&\qquad\qquad\times\frac{1}{\sqrt{2}}\frac{\Gamma(2l+4n+5/2)}
{(2\pi)^{2l+4n+5/2}}\\
&=\frac{\sqrt{|\bs{b}_S|}}{2\sqrt{2}\cdot\pi^{\frac{3}{2}}}\cdot
\frac{1}{\pi^{2l}}\cdot\frac{(l+2n+1/4)(l+2n+3/4)}{(l+2n+1/2)^2}
\cdot\Gamma(l+1)\Gamma(l+2)\\
&\qquad\qquad\times\frac{\Gamma(l+2n+1/4)\Gamma(l+2n+3/4)}
{\Gamma(l+2n+1)\Gamma(l+2n+1/2)}
\end{align*}
\end{proof}

{}From the above arguments we have
\begin{prop}
\begin{align*}
&\tilde{T}\circ T^{-1} = \frac{\bs{c}_l}{\bs{a}_l}Id\\
&=\sqrt{\frac{|\bs{b}_S|}{|\bs{b}_{\Ha}|}}\frac{1}{2^{\frac{5}{4}}}\cdot
\frac{(l+2n+1/4)(l+2n+3/4)}{(l+2n)(l+2n+1/2)}\\
&\qquad\qquad\times\frac{\Gm(l+2n+1/4)\Gm(l+2n+3/4)}{\Gm(l+2n+1/2)^2}Id
\end{align*}
on $H_l$.  Hence 
$\tilde{T}$ is an isomorphism between
$\mathfrak{h}^{\mathcal{G}}_{\Ha}$
and $L_2(P^n\Ha)$, because ${\bs{c}_l}/{\bs{a}_l}$ converges to
$\sqrt{\frac{|\bs{b}_S|}{|\bs{b}_{\Ha}|}}2^{-5/4}$ = $\pi /2$ 
when $l \to \infty$.
\end{prop}

Although the operator $\frac{\sqrt{\bs{a}_l}}{\bs{c}_l}\tilde{T}$ is
unitary on each $\mathcal{P}^{\Ha}_l$, no constant multiple of 
$\tilde{T}$ can be a unitary
operator between  $\mathfrak{h}^{\mathcal{G}}_{\Ha}$ and
$L_2(P^n\Ha)$.
However we have similar properties as for the operator $T$, that is
$\tilde{T}$ can also be understood as a quantization operator 
of the flow $\{\sigma_t\}$.

\begin{thm}
\begin{enumerate}
\item The operator $\tilde{T}$ is an isomorphism between the spaces
$\mathfrak{h}_{\Ha}^{\mathcal G}$ and $L_2(P^n\Ha)$.
\item The following diagram is commutative
\begin{equation}
\begin{CD}
\mathfrak{h}_{\Ha}^{\mathcal G} @>{\sigma_t^*}>>
\mathfrak{h}_{\Ha}^{\mathcal G}\\
@V{\tilde{T}}VV  @VV{\tilde{T}}V\\
L_2(P^n\Ha) @>>{e^{-t\sqrt{-1}\sqrt{\Delta^{\Ha}+(2n+1)^2}}}>
L_2(P^n\Ha).
\end{CD}
\end{equation}
\end{enumerate}
\end{thm}


\section{The reproducing kernel}

In this section we show that the Hilbert space
$\mathfrak{h}^{\mathcal{G}}_{\Ha}$
has a reproducing kernel.

We see easily that any function in
$\mathfrak{h}^{\mathcal{G}}_{\Ha}$
is holomorphic.  Let $f$ be in $\mathfrak{h}^{\mathcal{G}}_{\Ha}$ and
let $f = \sum_{l=0}^{\infty}f_l$ with $f_l \in \mathcal{P}^{\Ha}_l$,
then we have
\begin{equation*}
\mathcal{B}_l(f) = \mathcal{B}_l(f_l).
\end{equation*}
Hence we have
\begin{equation*}
f=\sum \frac{1}{\bs{b}_l}\mathcal{A}_l \circ \mathcal{B}_l(f).
\end{equation*}
Now we can rewrite this expression as
\begin{align*}
&f(A')\\
& = \sum_l \frac{1}{\bs{b}_l}
\int_{\widetilde{\mathbf{E}_{\Ha}}}f(A)\left(\int_{P^n\Ha}\Einner{P}{A'}_{\C}^l
\overline{\Einner{P}{A}_{\C}^l}\bs{v}_{\Ha}\right)\\
&\qquad\times e^{-2\sqrt[4]{2}\pi\sqrt{\norm{A}}}
\sqrt{|\bs{a}_{\Ha}|}\norm{A}^{n+1}\Omega_{\Ha}.
\end{align*}

Let 
\begin{equation}
\mathcal{R}(A,A')= \sum\frac{1}{\bs{b}_l}\int_{P^n\Ha}
\Einner{P}{A}_{\C}^l{\overline{\Einner{P}{A'}_{\C}^l}\bs{v}_{\Ha}},
\end{equation}
then we have
\begin{equation}
|\mathcal{R}(A,A')|\le \vol(P^n\Ha)\sum\frac{1}{2^l|\bs{b}_l|}\norm{A}^l\norm{A'}^l.
\end{equation}
{}From the expression of $\bs{b}_l$ (Proposition \ref{prop:3.5}) we
know that
\begin{equation}
\frac{\bs{b}_{l}}{\bs{b}_{l+1}} = O(l^{-4}).
\end{equation}
Hence $\mathcal{R}(A,\overline{A'})$ 
is a holomorphic function on
$\widetilde{\mathbf{E}_{\Ha}}\times
\widetilde{\mathbf{E}_{\Ha}}$. In fact it is holomorphic on all of the
space $M(2n+2,\C)$ and we have $\mathcal{R}(A,A')=\overline{\mathcal{R}(A',A)}$. 
By a similar estimation we have that for each fixed $A'\in M(2n+2,\C)$, 
the function 
$\mathcal{R}(A,A')$ of $A$ is in 
$\mathfrak{h}^{\mathcal{G}}_{\Ha}$.  In summary 
\begin{prop}
The Hilbert space $\mathfrak{h}^{\mathcal{G}}_{\Ha}$ has the
reproducing kernel
\begin{equation}
\mathcal{R}(A,A')= \sum\frac{1}{\bs{b}_l}\int_{P^n\Ha}
\Einner{P}{A}^l_{\C}{\overline{\Einner{P}{A'}^l_{\C}}\bs{v}_{\Ha}},
\end{equation}
and we have
\begin{equation*}
|f(A')|\le\norm{\mathcal{R}(A', \cdot)}\norm{f}.
\end{equation*}
\end{prop}

\bigskip

{\bf Acknowledgments}. {\it The author would like to express his
  hearty thanks to professor Maurice de Gosson and Mrs. Charlyne de
  Gosson for their kind hospitality during his stay at Blekinge
  Institute of Technology (Sweden), where the final version of this
  article was written. } 
\bigskip


\providecommand{\bysame}{\leavevmode\hbox to3em{\hrulefill}\thinspace}

\end{document}